\documentclass[12pt]{amsart}

\author{Paul \textsc{Poncet}}
\address{CMAP, \'{E}cole Polytechnique, Route de Saclay, 91128 Palaiseau Cedex, France \\
and INRIA, Saclay--\^{I}le-de-France}

\email{poncet@cmap.polytechnique.fr}

\usepackage{stmaryrd}
\usepackage{amssymb}
\usepackage{amsmath}
\usepackage{graphicx}
\usepackage{t1enc}
\usepackage[dvips]{color}
\usepackage[latin1]{inputenc}
\usepackage{yhmath}
\usepackage{setspace} 
\usepackage{fancybox} 
\usepackage{mathabx}
\usepackage{calrsfs} 
\newcommand{\call}{\mathrsfs}
\usepackage{times} 
\usepackage{ifpdf}
\def\twoheaddownarrow{\rlap{$\downarrow$}\raise-.5ex\hbox{$\downarrow$}}
\def\twoheaduparrow{\rlap{$\uparrow$}\raise.5ex\hbox{$\uparrow$}}

\newcommand{\Spec}{\mathsf{Spec}}

\DeclareMathOperator*{\im}{Im}

\newtheorem{theorem}{Theorem}[section]

\newtheorem{proposition}[theorem]{Proposition}
\newtheorem{lemma}[theorem]{Lemma}

\theoremstyle{definition}

\newenvironment{acknowledgements}[1][]{\par\vspace{0.5cm}\noindent\textbf{Acknowledgements#1.} }{\par}

\begin{document}
\setcounter{page}{1}

\title{
A class of compact subsets\\ for non-sober topological spaces}
\date{Decembre 30, 2009}

\subjclass[2010]{54D10, 
                 54D30, 
                 54D40, 
                 54D80, 
                 54C25} 

\keywords{topological space, compact saturated set, sober space, quasisober space, Hofmann--Mislove Theorem, capacities}

\begin{abstract}
We define a class of subsets of a topological space that coincides with the class of compact saturated subsets when the space is sober, and with enough good properties when the space is not sober. This class is introduced especially in view of applications to capacity theory. 
\end{abstract}

\maketitle
\tableofcontents


\section{Introduction}

It is well-known that in a sober (non-Hausdorff) space, compact subsets are advantageously replaced by the class of compact saturated subsets. The Hofmann--Mislove Theorem is an important result of general topology that illustrates this fact, for it gives a form of duality between Scott-open filters on the lattice of open subsets and compact saturated subsets in sober spaces. 
However, if the topological space is not sober, the result may fail. 
Since this duality is a very desirable property, especially for applications to capacity theory, a question arises: what collection could be used on every topological space to play the role of compact saturated subsets? 
This is the aim of this note to propose such a collection.

\section{Preliminaries on topological spaces}

Let $X$ be a topological space. The closure of a subset $A$ of $X$ is written $\overline{A}$, and for a singleton $\{x\}$ we write $\overline{x}$. 
The space $X$ satisfies the \textit{$T_0$ axiom} if $\overline{x} = \overline{y}$ implies $x = y$, for all $x, y \in X$. 
A subset $C$ of $X$ is \textit{irreducible} if it is nonempty and, for all closed subsets $F, F'$ of $X$, $C \subset F \cup F'$ implies $C \subset F$ or $C \subset F'$. Note that the closure of a point is always an irreducible closed set. 
We say that $X$ is \textit{quasisober} if every irreducible closed subset is the closure of a point, and \textit{sober} if it is quasisober and $T_0$. 
A subset of $X$ is \textit{saturated} if it is an intersection of open sets. 
In the general case ($X$ not necessarily $T_0$), consider the equivalence relation $x \sim y \Leftrightarrow \overline{x} = \overline{y}$, and let $X_0$ be the derived quotient space, endowed with the quotient topology, i.e.\ $G_0$ is open in $X_0$ if $\pi_0^{-1}(G_0)$ is open in $X$, where $\pi_0$ be the quotient map $x \mapsto \pi_0(x) := \{ y \in X : \overline{x} = \overline{y} \}$. For this topology, $X_0$ is $T_0$. 
Since $\pi_0$ is onto, we have $\pi_0(\pi_0^{-1}(A)) = A$ for all $A \subset X$. 
The map $A \mapsto \pi_0^{-1}(\pi_0(A))$ is a closure operator on $2^X$ called the \textit{$0$-closure}, and a closed set for this closure is \textit{$0$-closed}. Note that the class of $0$-closed subsets includes all closed subsets, all saturated subsets, hence in particular all open subsets. 

\section{Sobrification of topological spaces}

An important consequence of the Hofmann--Mislove Theorem evoked in the introduction is the following result, which emphasizes the r\^{o}le of compact saturated subsets for non-Hausdorff spaces. 

\begin{theorem}\label{dual}
Assume that $X$ is a sober topological space, and let $(Q_j)_{j \in J}$ be a filtered family of compact saturated subsets. Then the following statements hold:
\begin{enumerate}
  \item the collection of compact saturated subsets is closed under finite unions, 
	\item\label{dual1} $\bigcap_{j \in J} Q_j$ is compact saturated, 
	\item\label{dual2} if $\bigcap_{j \in J} Q_j \subset G$ for some open $G$, then $Q_j \subset G$ for some $j \in J$. 
\end{enumerate}
\end{theorem}


Note that (\ref{dual2}) implies (\ref{dual1}). 
This theorem was proved by Hofmann and Mislove \cite{Hofmann81}. A different proof is due to Keimel and Paseka \cite{Keimel94}. See Kov\'ar \cite{Kovar04} for an extension to generalized topological spaces, and Jung and S\"{u}nderhauf \cite{Jung96} for an enlightenment of this result in the context of proximity lattices. Also, Norberg and Vervaat \cite{Norberg97b} successfully applied this result, in a non-Hausdorff setting, to the theory of capacities which dates back to Choquet \cite{Choquet54}. 

As remarked in \cite{Norberg97b}, if $X$ is not sober, (\ref{dual2}) may fail. 
Before proposing a slight modification of the collection of compact saturated subsets for non-sober spaces, we need a result due to Hofmann and Lawson \cite{Hofmann78}, with a few additional features. We denote by $\call{O}(X)$ (resp.\ $\call{F}(X)$) the class of open (resp.\ closed) subsets of the space $X$. 

\begin{proposition}\label{psi}
Let $X$ be a topological space. There is a sober space $S$ and a continuous map $\xi : X \rightarrow S$ such that 
\begin{enumerate}
	\item\label{psi0} $\xi$ is open onto its image. 
	\item\label{psi00} $\xi$ is closed onto its image. 
	\item\label{psi0bis} There is some lattice isomorphism $\sigma : \call{O}(X) \rightarrow \call{O}(S)$ such that $\xi(G) = \sigma(G) \cap \im\xi$ and $G = \xi^{-1}(\sigma(G))$ for all $G\in \call{O}(X)$. 
	\item\label{psi00bis} There is some lattice isomorphism $\tau : \call{F}(X) \rightarrow \call{F}(S)$ such that $\xi(F) = \tau(F) \cap \im\xi$ and $F = \xi^{-1}(\tau(F))$ for all $F\in \call{F}(X)$. 	
	\item\label{psi1} $A = \xi^{-1}(\xi(A))$, for all $0$-closed subsets $A$ of $X$. 
	\item\label{psi2} $\xi(\bigcap_{j \in J} A_j) = \bigcap_{j \in J} \xi(A_j)$, for all families $(A_j)_{j \in J}$ of $0$-closed subsets of $X$. 
	\item\label{psi3} $\xi$ is injective if and only if $\xi$ is an embedding if and only if $X$ is $T_0$. 
	\item\label{psi4} $\xi$ is surjective if and only if $\xi$ is a quotient map if and only if $X$ is quasisober. 
\end{enumerate}
\end{proposition}

\begin{proof}
Let $S = \Spec \call{O}(X)$ be the sobrification of $X$, which consists of the set of prime elements of $\call{O}(X)$ ordered by inclusion. We equip $S$ with the hull-kernel topology. See \cite{Hofmann78}, \cite{Hofmann81}, and the monograph by Gierz et al.\ \cite{Gierz03}, for more details. 
Choose $\xi$ as the sobrification map $X \ni x \rightarrow X \setminus \overline{x} \in S$. 
Then \cite[Proposition~2.7]{Hofmann78} tells us that $\xi$ is continuous, open onto its image, and that (\ref{psi0bis}) and (\ref{psi3}) are satisfied. 

(\ref{psi00}) and (\ref{psi00bis}) Define $\tau$ by $\tau(F) = \{ P \in S : X \setminus F \subset P \} = S \setminus \sigma(X \setminus F)$ for every closed subset $F$ of $X$. Then $\xi(x) \in \tau(F) \cap \im\xi$ iff $X \setminus F \subset X\setminus\overline{x}$ iff $x \in F$ iff $\xi(x) \in \xi(F)$. (The fact that $\xi(x) \in \xi(F)$ implies $x \in F$ can be seen as follows. If $\xi(x) \in \xi(F)$, there is some $f \in F$ such that $\xi(x) = \xi(f)$, thus $\overline{x} = \overline{f} \subset F$, hence $x \in F$). This proves that $\tau(F) \cap \im\xi = \xi(F)$ for every closed subset $F$ of $X$. This also shows that $\xi$ is closed onto its image. The equality $F = \xi^{-1}(\tau(F))$ is a direct consequence of the definitions. 

(\ref{psi1}) Let $A \subset X$ be $0$-closed. 
We show that $A \supset \xi^{-1}(\xi(A))$ (the reverse inclusion is always true). Let $x \in \xi^{-1}(\xi(A))$. There is some $y \in A$ such that $\xi(x) = \xi(y)$, i.e.\ $X\setminus \overline{x} = X \setminus \overline{y}$. This implies $x \sim y$, hence $x \in \subset A$ since $A$ is $0$-closed. 

(\ref{psi2}) Let $(A_j)_{j \in J}$ be some family of $0$-closed subsets of $X$, and let us prove the inclusion $\xi(\bigcap_{j \in J} A_j) \supset \bigcap_{j \in J} \xi(A_j)$ (the reverse inclusion is trivially satisfied). Let $y \in \bigcap_{j \in J} \xi(A_j)$, and let $j_0 \in J$. For all $j \in J$, there is some $x_j \in A_j$ such that $y = \xi(x_j)$. Since $\xi(x_{j_0}) = \xi(x_j)$, $x_{j_0} \in \xi^{-1}(\xi(A_j)) = A_j$, hence $x_{j_0} \in \bigcap_{j \in J} A_j$, and $y \in \xi(\bigcap_{j \in J} A_j)$. 

(\ref{psi4}) The map $\xi$ is surjective if and only if every prime element of $\call{O}(X)$ is of the form $X\setminus\overline{x}$ for some $x \in X$, if and only if every irreducible closed subset of $X$ is of the form $\overline{x}$ for some $x \in X$, which exactly means that $X$ is quasisober. Moreover if $\xi$ is surjective, then $\xi$ is open, hence a quotient map.  
\end{proof}

The class we are seeking for is given by the next theorem. 

\begin{theorem}\label{thmsuperdualbis}
Let $X$ be a topological space, and define the class 
$$
\call{R}(X) = \{ R \subset X : R \mbox{ is $0$-closed in $X$ and } \xi(R) \mbox{ is compact saturated in $S$ } \}. 
$$
Then every element of $\call{R}(X)$ is compact saturated in $X$, and if $(R_j)_{j \in J}$ is a filtered family of elements of $\call{R}(X)$, then the  following statements hold: 
\begin{enumerate}
	\item\label{thmsuperdualbis2} $\call{R}(X)$ is closed under finite unions, 
	\item\label{thmsuperdualbis3} $\bigcap_{j \in J} R_j \in \call{R}(X)$, 
	\item\label{thmsuperdualbis4} if $\bigcap_{j \in J} R_j \subset G$ for some open $G$, then $R_j \subset G$ for some $j \in J$. 
\end{enumerate}
\end{theorem}

\begin{proof}
First, let $R$ be a $0$-closed subset of $X$ with $\xi(R)$ compact saturated in $S$. To show that $R$ is saturated in $X$, write 
$$
R = \xi^{-1}(\xi(R)) = \xi^{-1}(\bigcap_{U \supset \xi(R)} U) =  \bigcap_{U \supset \xi(R)} \xi^{-1}(U), 
$$
which is saturated as an intersection of open sets. 
Now we show that $R$ is compact, so assume that $R \subset \bigcup_{j\in J} G_j$, for open subsets $G_j$ of $X$. Then $\xi(R) \subset \xi(\bigcup_{j\in J} G_j) \subset \bigcup_{j\in J} \sigma(G_j)$. Since $\sigma(G_j)$ is open and $\xi(R)$ is compact in $S$, we deduce $\xi(R) \subset \bigcup_{j\in K} \sigma(G_j)$, for some finite $K \subset J$. Thus $R \subset \xi^{-1}(\bigcup_{j\in K} \sigma(G_j)) = \bigcup_{j\in K} \xi^{-1}(\sigma(G_j)) = \bigcup_{j\in K} G_j$, and $R \subset \bigcup_{j\in K} G_j$, so we conclude that $R$ is compact saturated in $X$. 

(\ref{thmsuperdualbis2}) is obvious, and (\ref{thmsuperdualbis3}) is implied by (\ref{thmsuperdualbis4}). 

(\ref{thmsuperdualbis4}) Assume that $\bigcap_{j \in J} R_j \subset G$ for some open $G \subset X$. We get 
$$
\bigcap_{j \in J} \xi(R_j) = \xi(\bigcap_{j \in J} R_j) \subset \sigma(G)
$$ 
thanks to Proposition~\ref{psi}. Since $\xi(R_j)$ is compact saturated in $S$ and $\sigma(G)$ is open in $S$, we apply Theorem~\ref{dual} to the sober space $S$, and we have $\xi(R_{j_0}) \subset \sigma(G)$ for some $j_0 \in J$. We get $R_{j_0} \subset \xi^{-1}(\sigma(G)) = G$, and (\ref{thmsuperdualbis4}) is proved. 
\end{proof}

We conclude this note by showing that the class introduced coincides with the class of compact saturated subsets if and only if $X$ is quasisober.

\begin{lemma}\label{lem:quasisob}
Let $S$ be a quasisober topological space. Then every saturated subset of $S$ is quasisober. 
\end{lemma}

\begin{proof}
Let $A$ be a saturated subset of $S$, and let $C$ be a nonempty irreducible closed subset of $A$. We can write $C = A \cap \overline{C}$, where $\overline{C}$ is the closure of $C$ in $S$, and it is easy to see that $\overline{C}$ is irreducible in $S$. Hence there is some $x \in S$ such that $\overline{C} = \overline{x}$, and $C = A \cap \overline{x}$. The proof will be over if we show that $x \in A$. So assume that $x \notin A$. Since $A$ is saturated, there exists some open subset $G$ of $S$ containing $A$ with $x \notin G$, i.e.\ $\overline{x} \subset S \setminus G$. Hence, $C \subset A \cap (S\setminus G) = \emptyset$, a contradiction. 
\end{proof}

\begin{theorem}
Let $X$ be a topological space. Then the following conditions are equivalent:
\begin{enumerate}
	\item\label{thm:fin1} $X$ is quasisober, 
	\item\label{thm:fin2} $\call{R}(X)$ covers $X$, 
	\item\label{thm:fin3} $\call{R}(X)$ coincides with the class of compact saturated subsets of $X$. 
\end{enumerate}
\end{theorem}


\begin{proof}
We write $\call{Q}(X)$ for the collection of compact saturated subsets of $X$. At first we show that (\ref{thm:fin3}) $\Rightarrow$ (\ref{thm:fin2}) $\Rightarrow (\im\xi \mbox{ saturated in } S) \Rightarrow$ (\ref{thm:fin3}), where $S$ and $\xi$ are defined by Proposition~\ref{psi}.  
If $\call{R}(X) = \call{Q}(X)$, then $X$ is obviously covered by $\call{R}(X)$. 
And if $X$ is covered by $\call{R}(X)$, then $\im\xi$ is equal to $\bigcup_{R} \xi(R)$, hence is saturated in $S$ as a union of saturated subsets of $S$. 
Now assume that $\im\xi$ is saturated in $S$, and let us show that $\call{Q}(X) = \call{R}(X)$. So let $Q \in \call{Q}(X)$. Then $Q$ is $0$-closed like every saturated subset, and since $\xi$ is continuous, $\xi(Q)$ is compact in $S$. So it remains show that $\xi(Q)$ is saturated in $S$. This is indeed the case, for $\xi(Q) = \xi(\bigcap_{G \supset Q} G) =  \bigcap_{G \supset Q} \sigma(G) \cap \im\xi$ is saturated as an intersection of saturated subsets. 

To finish the proof, we show that (\ref{thm:fin1}) $\Leftrightarrow (\im\xi \mbox{ saturated in } S)$. If $X$ is quasisober, then $\xi$ is surjective, hence $\im\xi = S$ is saturated in $S$. 
Conversely, if $\im\xi$ is saturated in $S$, then it is quasisober by Lemma~\ref{lem:quasisob}. Let us show that $X$ is quasisober. For this purpose, pick some nonempty closed irreducible subset $C$ of $X$. Then $\xi(C)$ is clearly an irreducible subset of $\im\xi$, and it is also closed in $\im\xi$ since $\xi$ is closed onto its image. Hence we can write $\xi(C) = \overline{z} \cap \im\xi$ for some $z \in \im\xi$, i.e.\ $\xi(C) = \overline{\xi(x)} \cap \im\xi$ for some $x \in X$. By continuity of $\xi$, we have $\xi(C) = \xi(\overline{x})$. But the closed subsets $C$ and $\overline{x}$ are $0$-closed, hence $C = \xi^{-1}(\xi(C)) = \xi^{-1}(\xi(\overline{x})) = \overline{x}$, and we have proved that $X$ is quasisober. 
\end{proof}




\begin{acknowledgements}
I wish to thank N.\ Noble for pointing out a mistake in an earlier version of the proof of Lemma~\ref{lem:quasisob}. 
\end{acknowledgements}

\bibliographystyle{plain}

\def\cprime{$'$} \def\cprime{$'$} \def\cprime{$'$} \def\cprime{$'$}
  \def\ocirc#1{\ifmmode\setbox0=\hbox{$#1$}\dimen0=\ht0 \advance\dimen0
  by1pt\rlap{\hbox to\wd0{\hss\raise\dimen0
  \hbox{\hskip.2em$\scriptscriptstyle\circ$}\hss}}#1\else {\accent"17 #1}\fi}
  \def\ocirc#1{\ifmmode\setbox0=\hbox{$#1$}\dimen0=\ht0 \advance\dimen0
  by1pt\rlap{\hbox to\wd0{\hss\raise\dimen0
  \hbox{\hskip.2em$\scriptscriptstyle\circ$}\hss}}#1\else {\accent"17 #1}\fi}

\end{document}